\documentclass[12pt]{amsart}

\usepackage{amssymb}

\newtheorem{thm}{Theorem}%[section]

\newtheorem{lem}[thm]{Lemma}
\newtheorem{cor}[thm]{Corollary}

\newtheorem{prop}[thm]{Proposition}

\newtheorem{conj}[thm]{Conjecture}
\theoremstyle{definition}
\newtheorem{defn}[thm]{Definition}

\newtheorem{say}[thm]{}
\newtheorem{exmp}[thm]{Example}

\newtheorem{ques}[thm]{Question}    %!!!!!!!!!!!!!!!!!!!!

\newtheorem{ack}{Acknowledgments}

\newtheorem{defn-thm}[thm]{Definition--Theorem}  %!!!!!!!!!!!!!!!!!!!!!!!!
\theoremstyle{remark}

\renewcommand{\a}[0]{{\mathbb A}}

\newcommand{\p}[0]{{\mathbb P}}
\newcommand{\f}[0]{{\mathbb F}}
\newcommand{\q}[0]{{\mathbb Q}}
\newcommand{\map}[0]{\dasharrow}
\newcommand{\qtq}[1]{\quad\mbox{#1}\quad}

\newcommand{\chr}[0]{\operatorname{char}}

\begin{document}
\bibliographystyle{amsplain}

\title{Unirationality of cubic hypersurfaces}
\author{J\'anos Koll\'ar}

\maketitle

A remarkable result of \cite{segre} says that a 
smooth cubic surface
over $\q$ is unirational iff it has  a rational point. 
\cite[II.2]{manin} observed that similar arguments work
for higher dimensional cubic hypersurfaces satisfying
a certain genericity assumption over any infinite  field.
\cite[2.3.1]{ctssd} extended the result of Segre to  any normal
cubic hypersurface (other than cones) over a field of characteristic zero.
It is also clear that  the result should
hold for all
sufficiently large finite fields,
though the details were not worked out in general.
\cite[IV.8]{manin}
settles the cubic surface case for  finite fields with at least
34 elements. The aim of this note is to observe that
a variant of the Segre--Manin method works
for all  fields and for all cubics:

\begin{thm}\label{smooth-main.thm} Let $k$ be a  field and 
$X\subset \p^{n+1}$  a smooth  cubic hypersurface of
dimension $n\geq 2$
over $k$. 
Then the following are equivalent:
\begin{enumerate}
\item $X$ is unirational (over $k$).
\item $X$ has  a $k$-point.
\end{enumerate}
\end{thm}

Similar results hold for singular cubic hypersurfaces,
with a few exceptions.

\begin{thm}\label{main.thm} Let $k$ be a perfect field and 
$X\subset \p^{n+1}$  an irreducible  cubic hypersurface of
dimension $n\geq 2$
over $k$ which is not a cone over an $(n-1)$-dimensional cubic. 
Then the following are equivalent:
\begin{enumerate}
\item $X$ is unirational (over $k$).
\item $X$ has  a $k$-point.
\item $X$ has a smooth $k$-point.
\end{enumerate}
\end{thm}

\begin{say}[Nonperfect fields] {\ }
 Over  nonperfect fields of characteristic
3, there are nonsingular cubic hypersurfaces of arbitrary
dimension  which are not unirational but do have a $k$ point
 (\ref{char3.exmp}).  I have not been able to find examples
 in characteristic 2.
\end{say}

\begin{ques}
Unirationality of varieties is very poorly understood in general
and there are very basic open questions. We do not
even have a list of  unirational surfaces and
very few examples are known in higher dimensions. For instance, 
let $X$ be a smooth projective variety over $k$
such that $X$ is unirational over $\bar k$.
 Assume for simplicity that $k$ is  infinite and
consider the following properties:
\begin{enumerate}
\item $X$ is unirational (over $k$).
\item $X(k)$ is dense in $X$.
\item $X$ has  a $k$-point.
\end{enumerate}
It is clear that each property implies the next. They are
equivalent for cubic hypersurfaces by (\ref{smooth-main.thm}).
It is extremely unlikely that they are always equivalent,
but no counter examples are known.
\end{ques}

\begin{say}[Proof of (\ref{main.thm}.1) $\Rightarrow$
 (\ref{main.thm}.2) $\Rightarrow$ (\ref{main.thm}.3)]{\ }

 (\ref{main.thm}.1) $\Rightarrow$ (\ref{main.thm}.2)
 is clear for infinite fields.
For finite fields it 
follows from \cite{nishimura}.

Assume (\ref{main.thm}.2) and let $x\in X$ be a $k$-point. We are
 done if $X$ is smooth at $x$. Otherwise $x$ is a double point and we can 
 choose affine coordinates such that $x=(0,\dots,0)$ and $X$ is given 
by an equation
$ q(x_1,\dots,x_n)+c(x_1,\dots,x_n)=0$ where $q$ is quadratic and $c$ is
 cubic. Assume that there is a point $(p_1,\dots,p_n)\in k^n$
such that $q(p_1,\dots,p_n)\neq 0$. 
Then the line connecting the origin and $(p_1,\dots,p_n)$
intersects $X$ in a single point outside the origin and this is
a smooth $k$-point of $X$. Thus we are done unless
$q$ vanishes everywhere
on $k^n$.

However, if  a homogeneous polynomial $f$ of degree $d$  vanishes 
on $k^n$ and $|k|\geq d$ then $f$ is identically zero. \qed
\end{say}

The intersting  part is to show unirationality
starting with a smooth $k$-point. The construction is
presented in 3 stages, successive  version working in
greater and greater generality. At least in retrospect,
all of this is only a slight modification of the works of Segre.

\begin{say}[First unirationality construction]{\ }

Let $X\subset \p^{n+1}$ be a cubic and $p\in X$ a point.
Let $C_p$  denote the intersection of
$X$ with the tangent plane at $p$.
We expect that usually  $C_p$
is an irreducibe cubic with a double point at $p$.
If this is indeed the case then the inverse of the projection from $p$
gives a birational map $ \pi_p: \p^{n-1}\map C_p$. If $p\in X(k)$
then $C_p$ is birational to $\p^{n-1}$ over $k$.

Assume next that we have   two points $p,q\in X$
and  $C_p,C_q$ 
are  both irreducible  with a double point at $p$ (resp.\ $q$).
Define the ``3rd intersection point'' map 
$$
\phi:C_p\times C_q\map X
$$
as follows. Take $u\in C_p, v\in C_q$. If the line connecting
$u,v$  is not contained in $X$, it has a unique 3rd intersection
 point with $X$; call it $\phi(u,v)$. Under very mild
genericity assumptions (\ref{dom.map.lem}) this is a well defined dominant
map. Thus we get that $X$ is unirational via
$$
\Phi: \p^{n-1}\times \p^{n-1}
\stackrel{\pi_p\times \pi_q}{\map} 
C_p\times C_q\stackrel{\phi}{\map} X.
$$
\end{say}

\begin{defn}[Restriction of scalars]{\ }

Let $L/K$ be a finite degree field extension. 
Restriction of scalars (or Weil restriction) is  way to associate
to an $L$-variety $U$  a $K$-variety ${\frak R}_{L/K}U$ 
such that there is a natural identification
of the $L$-points of $U$ with the
$K$-points of ${\frak R}_{L/K}U$. 
This dictates that  $\dim {\frak R}_{L/K}U=\deg(L/K)\cdot \dim U$;
 see \cite[7.6]{blr} for details.

This can be done very explicitly in the affine case as follows.
Let $U\subset \a^n$ be an affine variety defined over
$L$.   
Choose equations $f_i(x_1,\dots,x_n)$ for $U$
and let $e_1,\dots,e_d\in L$ be a $K$-basis. Choose new coordinates
$y_{ij}:i=1,\dots,n, j=1,\dots, d$ and set
$x_i=\sum_j e_jy_{ij}$. We can then write
$$
f_k(x_1,\dots,x_n)=\sum_{\ell} e_{\ell}f_{k\ell}(y_{ij})
\qtq{where $f_{k\ell}\in K[y_{ij}]$.}
$$
Let 
${\frak R}_{L/K}U$ be the subvariety of $\a^{nd}$ defined
by the equations $f_{k\ell}=0$.

In particular we see that ${\frak R}_{L/K}(\a^n)\cong \a^{dn}$.
In the projective case, ${\frak R}_{L/K}(\p^n)$ is
inconvenient  to describe by explicit equations but we at least
get that  ${\frak R}_{L/K}(\p^n)$ is birational to
$\p^{dn}$ over $K$. (They are not isomorphic for $d>1$.)
\end{defn}

\begin{say}[Second unirationality construction]{\ }

Let $X\subset \p^{n+1}$ be a cubic 
defined over $k$  and $k'\supset k$ a quadratic extension.
Let  $p\in X(k')$ be a point and $\bar p\in X(k')$ its conjugate.
(Let us ignore  that $k'\supset k$ may be inseparable
in characteristic 2.) We have conjugate birational maps
$\pi_p:\p^{n-1}\map C_p$ and $\pi_{\bar p}:\p^{n-1}\map C_{\bar p}$.
If $u\in \p^{n-1}(k')$ then $\pi_p(u)$ and
$\pi_{\bar p}(\bar u)$ are conjugate points of $X$,
thus the line connecting them is defined over $k$.
Hence $\phi(u,\bar u)\in X(k)$.
Putting this invariantly,  we obtain a rational map (defined over $k$)
$$
\Phi:{\frak R}_{k'/k}\p^{n-1} \map X
$$
which is dominant under mild genericity assumptions.
\end{say}

\begin{say}[Final unirationality construction]{\ }

Assume now that $X$ is a cubic defined over $k$
and $x\in X$ is a smooth $k$-point. Let $L$ be a line through
$x$. If $L$ is not contained in $X$ then it intersects
$X$ in a  point pair $\{p,q\}$. These points are usually not in $k$,
but they are conjugate over $k$ and lie in a quadratic extension $k'=k'(L)$
of $k$.
Hence, under some genericity assumptions,
we obtain a dominant map
$$
\Phi:{\frak R}_{k'/k}\p^{n-1} \map X
$$
which shows that $X$ is unirational. There are very few problems
if $k$ is infinite, since then a general choice of $L$
should work. (Extra work is needed in characteristic 2.)
The situation is less clear  over finite fields since there
may not be enough room to choose $L$ general; see, for example, 
(\ref{fewpts.exmp}).

To avoid this difficulty, we do not choose any line, rather
we work with all lines simultaneously. We should obtain a
 map
$$
\Psi:\cup_{x\in L\subset \p^{n+1}} {\frak R}_{k'(L)/k}\p^{n-1} \map X.
$$
We are in good shape if we can  indentify the left hand side
with a product $\p^n\times \p^{n-1}\times \p^{n-1}$,
at least birationally.
 Once this problem is settled,
it is enough to check dominance over the algebraic closure
where the previous arguments work. It seems best to give an
explicit algebraic description.
\end{say}

\begin{say}[Algebraic description of $\Psi$]{\ }

 We work in affine coordinates, assuming that
the origin is a smooth point of  $X$. Thus the equation of $X$
can be written as
$$
F=L(x_1,\dots,x_{n+1})+Q(x_1,\dots,x_{n+1})+C(x_1,\dots,x_{n+1})
$$
where $L$ is linear, $Q$ is quadratic and $C$ is cubic.
We may assume that $\partial F/\partial x_{n+1}$
is not indetically zero (for instance we can even assume that
 $L= x_{n+1}$).

We write down a  rational map
$$
\Psi: \a^{3n-2}(u_1,\dots,u_{n},v_1,\dots,v_{n-1},w_1,\dots,w_{n-1})
\map X.
$$
Later we check that it is dominant with a few exceptions.

Consider the universal line through the origin
$(\tau u_1,\dots, \tau u_{n}, \tau)$. It intersects $X$ in
two further points which correspond to the roots
of the quadratic equation
$$
L(u_1,\dots,u_{n},1)+\tau Q(u_1,\dots,u_{n},1)+
\tau^2C(u_1,\dots,u_{n},1)=0.
$$
The equation is irreducible if $X$ is irreducible.
Let its roots  be $t_1,t_2\in \overline{k(u_1,\dots,u_{n})}$. 

 The equation of the tangent space
of $X$ at ${\mathbf p}=(p_1,\dots,p_{n+1})\in X$ is
$$
\frac{\partial F}{\partial x_1}({\mathbf p})(x_1-p_1)+\cdots +
\frac{\partial F}{\partial x_{n+1}}({\mathbf p})(x_{n+1}-p_{n+1})=0.
$$
Thus the universal tangent line at 
$(t_1 u_1,\dots, t_1 u_{n}, t_1)$
 can be described parametrically as
$$
\begin{array}{l}
 x_1=t_1u_1+\sigma(v_1+t_1w_1),\dots,
x_{n-1}=t_1u_{n-1}+\sigma(v_{n-1}+t_1w_{n-1})\\
x_{n}=t_1u_{n}+\sigma,\\
 x_{n+1}=t_1-\sigma
\bigl(
\frac{\partial F}{\partial x_{n+1}}(t_1{\mathbf u},t_1)
\bigr)^{-1}
\sum_{i=1}^{n}
\frac{\partial F}{\partial x_i}(t_1{\mathbf u},t_1)(v_i+t_1w_i),
\end{array}
$$
where we set $v_n=1, w_n=0$.
Substituting the above parametric representation into $F$,
we obtain a cubic equation in $\sigma$
$$
\sum_{j=0}^3 \sigma^jH_j\qtq{where}
H_j\in k({\mathbf u},{\mathbf v},{\mathbf w}, t_1).
$$
 $H_0=H_1=0$ since  we have a tangent line, thus
the 3rd intersection point corresponds to the value
$\sigma=-H_2/H_3$. Thus we obtain a point
$$
Q_1\in k({\mathbf u},{\mathbf v},{\mathbf w}, t_1)^{n+1}.
$$
Replacing $t_1$ by its conjugate $t_2$ we obtain another point
$Q_2$.
The line connecting $Q_1$ and $Q_2$ can be given parametrically as
$$
L(\lambda)=\frac{\lambda-t_2}{t_1-t_2}Q_1+ 
\frac{\lambda-t_1}{t_2-t_1}Q_2,
$$
and this is  a parametrization over
 $k({\mathbf u},{\mathbf v},{\mathbf w})$.
Evaluating $F$ on the line we have that $t_1,t_2$ are roots,
so
$$
F(L(\lambda))=(A\lambda+B)
(C\lambda^2+Q\lambda+L).
$$
Thus if we expand
$$
F(L(\lambda))=
\sum_{j=0}^3 \lambda^jG_j,\qtq{then}
G_j\in k({\mathbf u},{\mathbf v},{\mathbf w}),
$$
and the 3rd root is
$$
-\frac{B}{A}=-\frac{G_2}{G_3}+\frac{Q}{C}.
$$
Substituting this into the parametrization of the line
gives
$$
\Psi({\mathbf u},{\mathbf v},{\mathbf w})
\in X(k({\mathbf u},{\mathbf v},{\mathbf w})).
$$
\end{say}

Depending on our definition of unirationality, we also
need to check the following:

\begin{lem} For a $k$-variety $X$ the following are equivalent:
\begin{enumerate}
\item There is a dominant map $\phi_m:\a^m\map X$ for some $m$.
\item There is a dominant map $\phi_m:\a^m\map X$ for $m=\dim X$.
\end{enumerate}
\end{lem}

Proof. Assume that $m>\dim X$.
There is a dense open set $U\subset \a^m$
such that $\phi_m|_U$ is open with $m-\dim X$ dimensional fibers.
Let ${\mathbf u}\in U$ be a point. 
If ${\mathbf u}\in Z\subset \a^m$ is a hypersurface  which
does not contain the irreducible component of the fiber of $\phi_m|_U$ 
through ${\mathbf u}$,
then $\phi_m|_Z:Z\map X$ is dominant.

If $k$ is infinite and $m>\dim X$ then
we can choose $Z$ to be a general 
 hyperplane.

Assume next that $k$ is finite. Fix a prime $\ell\neq \chr k$
and let $k'$ be the composite of all algebraic extensions of
degree $\ell^s$ of $k$. $k'$ is infinite, hence we can choose
a  point ${\mathbf u}=(u_1,\dots,u_m)\in U(k')$.
By permuting the coordinates we may assume that
$\deg k(u_m)/k\leq \deg k(u_1)/k$, or, 
equivalently, $k(u_m)\subset k(u_1)$.
This implies that
$u_m$ can be written as  a polynomial of $u_1$, hence
the ideal $I({\mathbf u})\subset k[x_1,\dots,x_m]$ contains an polynomial
of the form $x_m-p(x_1)$. This implies that
$I({\mathbf u})$ is generated by polynomials 
of the form $x_m-P(x_1,\dots,x_{m-1})$. Thus we can choose 
$Z=(x_m=P(x_1,\dots,x_{m-1}))$ for suitable $P$.\qed

\begin{say}[Proof of (\ref{main.thm}.3) $\Rightarrow$
 (\ref{main.thm}.1)]{\ }

Let $X$ be an irreducible  cubic hypersurface.
The set of all triple points of $X(\bar k)$ is a linear space
and it is defined over $k$ if $k$ is perfect. Thus 
$X$ is a cone over a cubic hypersurface without triple points.
Therefore, $X$ has no  triple points over $\bar k$.

Assume next that $X$ is not normal. The nonnormal locus has dimension 
$(n-1)$ and the linear space spanned by it is in $X$.
Thus the nonnormal locus is a linear space $L^{n-1}\subset \p^{n+1}$
which is defined over $k$ if $k$ is perfect.
Projecting form $L$ realizes $X$ as a $\p^{n-1}$-bundle over $\p^1$,
hence rational.

For the rest of the proof assume that $X$ is normal.
We need to check three conditions. 

First we prove that
  $C_x$ is irreducible with a double point at $x$
for general $x\in X(\bar k)$.
This is done in (\ref{good.C_p}).

Second, we need to check that the 3rd intersection point map
$\phi: C_p\times C_q\map X$ is dominant. 
It is, however, not enough to check this
for a general pair $p,q$. In our construction
$p,q$ are the two intersection points
of a line through $x$, hence dependent.
Assume that $\pi_x:X\map \p^n$, the projection from $x$, is separable.
Then for a generic line $x\in L$ we get 2 distinct intersection points
and both intersections are transverse. In particular,
the tangent space of $X$ at one point does not contain the other point.
In (\ref{dom.map.lem}) we see that 
this is sufficient to guarantee  that $\phi: C_p\times C_q\map X$ is dominant.

Third,  we need to consider the case when the projection
$\pi_x:X\map \p^n$ is inseparable. 
This can happen only in characteristic 2.
Over a perfect field a purely inseparable map induces
a purely inseparable map in the reverse direction,
hence in this case $X$ is (purely inseparably)
unirational. Nonetheless, we check in
(\ref{pi.cor}) that we can always choose a 
smooth $k$-point such that projection from it is separable.\qed
\end{say}

By looking at the proof we obtain the following for nonperfect
fields. We check in (\ref{proj.coord.lem}) that 
(\ref{nonperf.prop}.3) is satisfied for $X$ smooth. This
shows that our proof covers all smooth hypersurfaces.

\begin{prop}\label{nonperf.prop}
 Let $k$ be a nonperfect field and $X\subset \p^{n+1}$
a cubic hypersurface which is not a cone.
 Then the 3 parts of (\ref{main.thm}) are
equivalent if one of the following conditions  holds:
\begin{enumerate}
\item $\chr k\geq 5$.
\item $\chr k=3$ and $X$ has no triple points over $\bar k$.
\item $\chr k=2$ and there is a smooth $k$ point $p\in X$ such that 
projection from $p$ is
separable.\qed
\end{enumerate}
\end{prop}

\begin{prop}\label{good.C_p}
 Let $k$ be an algebraically closed field and
$X\subset \p^{n+1}$  a normal 
cubic hypersurface over $k$ without triple points.
Then $C_x$ is irreducible with a double point at $x$ for general $x\in X$.
\end{prop}

Proof.  Let $x\in X$ be arbitrary. If $C_x$ is irreducible with a
triple point at $x$ then $C_x$ is a cone, hence there is an 
$(n-2)$-dimensional family of lines through $x$. If $C_x$ is reducible
then either $C_x$ contains an $(n-1)$-dimensional linear space
 through $x$ or a quadric cone with vertex at $x$. In either case, there is an 
$(n-2)$-dimensional family of lines through $x$. Thus it is enough to prove
that for a general $x\in X$ the family of lines in $X$ through $x$
has dimension at most $n-3$. This is equivalent to proving
that a general surface section through $x$ has no lines through $x$.

$X$ has no triple points, hence by Bertini, a general surface section
of $X$ is also normal with
 no triple points.  

If $S$ is a normal cubic surface without
triple points then there are only finitely many lines through each
 double point. 
(Choose affine coordinates such that the equation becomes
$q(x_1,x_2,x_3)+c(x_1,x_2,x_3)=0$. The lines through
$(0,0,0)$ correspond to the solutions of $(q=c=0)\subset \p^2$.
If there are infinitely many solutions, then $q$ and $c$ have a 
common factor, thus the surface is reducible.)
Each line in the smooth locus has selfintersection
 $-1$, hence rigid.
Thus $S$ has only finitely many lines.\qed

\begin{lem}\label{dom.map.lem}
 Let $X\subset \p^{n+1}$ be an irreducible  cubic hypersurface.
Let $x,y\in X$ be smooth points and $C_x,C_y$ the corresponding
intersections with the tangent hyperplanes. Assume that
\begin{enumerate}
\item $C_x$ and $C_y$ are irreducible.
\item $x\not\in C_y$ and $y\not\in C_x$.
\end{enumerate}
Then the 3rd intersection point map $\phi:C_x\times C_y\map X$
is dominant.
\end{lem}

Proof.  Let us see first that $\phi$ is indeed defined.
Pick a point $u\in C_x$ which is a smooth point of $X$.
Pick $v\in C_y$ which is a smooth point of $X$ such that
$v$ does not lie on $T_uX$. If we now choose a general
$w\in C_x$ then $v$ does not lie on $T_wX$ and 
$w$ does not lie on $T_vX$. Thus the line connecting
$u,w$ has a unique third intersection point with $X$.
This shows that $\phi$ is defined at the pair $(v,w)$.

In order to prove dominance, 
we need to show that $\phi$ has at least one  fiber of dimension
$n-2$. Pick a point $z\in X$ which is not on
$C_x\cup C_y$ and 
let $\pi: \p^{n+1} \map T_yX$ denote the projection from $z$.
Then $\phi^{-1}(z)$ is the set of pairs $(v,w)$
such that $\pi(v)=w$. Thus 
we are done if 
$$
\dim ((C_y\cap \pi(C_x))\setminus (C_y\cap C_x))=n-2.
$$
For this it is sufficient to find one projection
$\pi': \p^{n+1} \map T_yX$ where this holds.
Then the same holds for a general projection and a general
projection always corresponds to a point of $X$.
 Pick any smooth point $v\in C_y$
and let $\pi'$ be a projection such that $\pi'(x)=v$. 
Then $\pi'(C_x)$ and $C_y$ intersect at $v$ but they have
 different multiplicty there. Hence their intersection has
 dimension $n-2$.\qed

\begin{exmp}\label{fewpts.exmp}
Let $S$ be the cubic surface $(x_0^3+x_1^3+x_2^3+x_3^3=0)$.
By \cite{hirschfeld} over the fields $\f_2,\f_4,\f_{16}$
all the points are on the 27 lines. Hence
the second unirationality construction  does not work over 
$\f_2$ and $\f_4$. (It does work over $\f_{16}$.)
\end{exmp}

\begin{exmp}\label{char3.exmp}
 Let $k$ be a field of characteristic 3 and
$t_i$ algebraically independent over $k$. Set $K=k(t_1,\dots, t_n)$ and 
$$
Y:=(y^3-yz^2=\sum_{i=1}^n t_ix_i^3)\subset \p^{n+1}.
$$
\begin{enumerate}
\item  $Y$ is non--singular.
\item Over $\bar K$, $Y$ is a cone over a cuspidal cubic curve.
\item $Y(K)=\{(0,1,0,\dots,0), (1,1,0,\dots,0), (1,-1,0,\dots,0)\}$.
\item $Y$ is not unirational (over $K$).
\end{enumerate}

Proof. $Y$ is the generic fiber of the
smooth variety 
$$
(y^3-yz^2=\sum_i t_ix_i^3)\subset
 \a^n_{(t_1,\dots,t_n)}\times \p^{n+1}_{(y,z,x_1,\dots,x_n)}
$$
over $\a^n$, thus $Y$ is non--singular. 
(2) holds since over $\bar K$ we can write our equation as
$$
(y-\sum_i \sqrt[3]{t_i}x_i)^3-yz^2=0.
$$

In order to see (3) 
we may as well assume that $k$ is algebraically closed.
Assume that we have relatively prime polynomials
$f,g,h_i\in k[t_1,\dots, t_n]$ such that
$$
f^3-fg^2=\sum_i t_ih_i^3.
$$
We are done if $h_1=\dots=h_n=0$. 
Otherwise, we can make  a substitution $t_i=c_it$ for $i=1,\dots, n$
and general $c_i$ to get a solution of
$$
f(f-g)(f+g)=t\cdot h^3 \qtq{with $f,g,h\in k[t]$ and $h\neq 0$.}
$$
We may assume that  $f$ and $g$ are relatively prime. 
Thus 2 of the factors $f,f-g,f+g$ are cubes and the third
is $t$ times a cube. However, 
$f+(f-g)+(f+g)=0$, hence if 2 are cubes
then so is their sum which is minus the 3rd factor.
 This is a contradiction.

Since $Y$ has only 3 points in $K$, it does not contain any
rational curves and so it is definitely not unirational.
\qed
\end{exmp}

\begin{lem}\label{proj.coord.lem}
  Assume that $\chr k=2$. Let $V\subset \p^{n+1}$
be the linear span of all points $p\in X(k)$ such that projection from
$p$ is a purely inseparable map $X\map \p^n$.
Let $(y_i=0)$ be equations of $V$ and $x_j$ coordinates on $V$.
Then the equation of $X$ can be written as
$$
f:=\sum_j \ell_j({\mathbf y})x_j^2+g({\mathbf y})
$$
where the $\ell_j$ are linear and $g$ is cubic.
If $V\neq\emptyset$ then $X$ is not smooth.
\end{lem}

Proof. We can choose coordinates such that
the points
$$
p_1=(1:0:\cdots:0), \dots, p_m=(0:\cdots :
\stackrel{m-th}{1}:0:\cdots :0)
$$
are in $X(k)$ and 
projection from
$p_i$ is a purely inseparable map for $i=1,\dots,m$. 
$p_i$ is inseparable iff $x_i$ occurs in the equation of
$X$ always with even exponent. This gives the above equation.

$\partial f/\partial x_j$ is zero, and the equations
$\partial f/\partial y_i=0$ have a common solution.
Since $f=3f=\sum_i (\partial f/\partial y_i)$, these give
singular points of $X(\bar k)$.
\qed

\begin{lem} Let $k$ be a perfect field of characteristic 2.
Let $X$ be a cubic of dimension at least 2 given by an equation
$$
f({\mathbf x},{\mathbf y}):=\sum_j \ell_j({\mathbf y})x_j^2+g({\mathbf y}).
$$
Then $X$ has a smooth $k$-point with nonzero $y$-coordinate.
\end{lem}

Proof. Assume first  that we have at least two $x$-variables.
If $\ell_1=c\ell_2$ then 
$$
\ell_1x_1^2+\ell_2x_2^2=\ell_1(x_1+\sqrt{c}x_2)^2
$$
thus we can change coordinates  to eliminate one $x$-variable.
Otherwise we can 
pick ${\mathbf y}_0$ such that $\ell_1({\mathbf y}_0)\neq 0$
and  $\ell_2({\mathbf y}_0)= 0$.
Then 
$$
(\sqrt{g({\mathbf y}_0)/\ell_1({\mathbf y}_0)},x_2,{\mathbf y}_0)
\in X(k)
$$
for every $x_2$.
Since 
$$
\frac{\partial f}{\partial y_i}=x_2^2+\frac{\partial g}{\partial y_i}
$$
the above point is smooth for suitable choice of $x_2$.

Thus assume that there is only one $x$-variable and
write the equation as $y_1x_1^2+g({\mathbf y})$.

Take any $(p_1,\dots,p_n)\in k^n$. If $p_1=0$
and $g(p_1,\dots,p_n)=0$
then $(x_1:p_1:\dots :p_n)\in X(k)$ for any $x_1$
and one of them is a smooth by looking at 
${\partial f}/{\partial y_1}$.

If $p_1\neq 0$ then $p_0:=\sqrt{-g(p_1,\dots,p_n)/p_1}\in k$
and  $(p_0:p_1:\dots :p_n)$ is a smooth point unless
$$
g(p_1,\dots,p_n)-p_1\frac{\partial g}{\partial y_1}(p_1,\dots,p_n)=0.
$$
Thus we are done unless the following holds:
\begin{enumerate}
\item $g-y_1(\partial g/\partial y_1)$ is nonzero for $y_1=0$, and
\item $g-y_1(\partial g/\partial y_1)$ is zero for $y_1\neq 0$.
\end{enumerate}
Write $g=\sum y_1^ig_{3-i}(y_2,\dots,y_n)$. Then
$$
g-y_1(\partial g/\partial y_1)=y_1^2g_{1}(y_2,\dots,y_n)+
g_{3}(y_2,\dots,y_n).
$$
$g_1$ is a linear form thus it has  a nontrivial zero
$(p_2,\dots,p_n)$.  If $g_3(p_2,\dots,p_n)=0$ then set $p_1=0$
and if $g_3(p_2,\dots,p_n)\neq 0$ then set $p_1=1$.\qed

\medskip

Combining the above lemmas we obtain:

\begin{cor} \label{pi.cor}
Let $k$ be a perfect field of characteristic 2 and
$X\subset \p^{n+1}$ a cubic with a smooth $k$-point. Assume that $n\geq 2$.
Then there is a smooth point
$x\in X(k)$ such that the projection from $x$ is separable.\qed
\end{cor}

\begin{ack}   I   thank J.-L.\ Colliot-Th\'el\`ene 
and J.\ Ellenberg
for helpful
comments and references.
Partial financial support was provided by  the NSF under grant number 
DMS-9970855. 
\end{ack}

\vskip1cm

\noindent Princeton University, Princeton NJ 08544-1000

\begin{verbatim}kollar@math.princeton.edu\end{verbatim}

\end{document}